\documentclass[12pt]{article}
\usepackage{verbatim}
\usepackage{fancyhdr}
\usepackage{graphicx}
\usepackage{mathrsfs}
\usepackage{amssymb}
\usepackage{pst-poly}
\usepackage{cite}

\newtheorem{lem}{Lemma}[section]
\newtheorem{thm}{Theorem}[section]

\newenvironment{pf}[1][Proof]{\noindent\textbf{#1.} }{\hfill\rule{1mm}{2mm}}

\makeatletter \@addtoreset{equation}{section} \makeatother

\begin{document}

\title{The bondage number of $(n-3)$-regular graphs of order $n$\thanks {The work was supported by NNSF
of China (No. 11071233).}}
\author
{Fu-Tao Hu, Jun-Ming Xu\footnote{Corresponding
author:
xujm@ustc.edu.cn}\\ \\
{\small Department of Mathematics}  \\
{\small University of Science and Technology of China}\\
{\small Hefei, Anhui, 230026, China} }
\date{}
 \maketitle

\begin{quotation}

\textbf{Abstract}: Let $G=(V,E)$ be a graph. A subset $D\subseteq V$
is a dominating set if every vertex not in $D$ is adjacent to a
vertex in $D$. The domination number of $G$ is the smallest
cardinality of a dominating set of $G$. The bondage number of a
nonempty graph $G$ is the smallest number of edges whose removal
from $G$ results in a graph with larger domination number of $G$. In
this paper, we determine that the exact value of the bondage number
of $(n-3)$-regular graph $G$ of order $n$ is $n-3$.

\vskip6pt\noindent{\bf Keywords}: domination number, bondage number, regular graph.

\noindent{\bf AMS Subject Classification: }\ 05C69

\end{quotation}

\section{Introduction}

For graph-theoretical terminology and notation not defined here, we
follow \cite{x03}. Specifically, a graph $G=(V,E)$ is considered as
an undirected graph without loops and multi-edges, where $V=V(G)$ is
the vertex-set and $E=E(G)$ is the edge-set. For a vertex $x$ in
$G$, let $N_G(x)=\{y\in V(G): xy\in E(G)\}$, $N_G[x]=N_G(x)\cup
\{x\}$ and $E_G(x)=\{xy: y\in N_G(x)\}$. The cardinality $|E_G(x)|$
is the degree of $x$. For two disjoint nonempty and proper subsets
$S$ and $T$ in $V(G)$, we use $E_G(S,T)$ to denote the set of edges
between $S$ and $T$ in $G$, and $G[S]$ to denote a subgraph of $G$
induced by $S$.

A vertex $y$ is said to be dominated by a vertex $x$ if $y\in
N_G(x)$ in graph $G$. A subset $D\subset V$ is a {\it dominating
set} of $G$ if $N_G(x)\cap D\ne \emptyset$ for every vertex $x$ in
$G-D$. The {\it domination number} of $G$, denoted by $\gamma(G)$,
is the minimum cardinality of all dominating sets of $G$.

In 1990, Fink {\it et al.}~\cite{fjkr90} introduced the bondage
number as a parameter for measuring the vulnerability of the
interconnection network under link failure. The {\it bondage number}
of a nonempty graph $G$, denoted by $b(G)$, is the minimum number of
edges whose removal from $G$ results in a graph with larger
domination number than $\gamma(G)$, that is,
 $$
 b(G) =\min\{|B| : B\subseteq E(G), \gamma(G - B)>\gamma(G)\}.
 $$
A nonempty subset $B\subseteq E(G)$ is said a {\it bondage set} of
$G$ if $\gamma (G-B)>\gamma(G)$. A bondage set $B$ is said to be
{\it minimum} if $|B|=b(G)$. In fact, if $B$ is a minimum bondage
set, then $\gamma(G-B)=\gamma(G)+1$, because the removal of one
single edge can not increase the domination number by more than one.

It is quite difficult to compute the exact value of the bondage
number for general graphs since it strongly depends on the
domination number of the graphs. Very recently, Hu and
Xu~\cite{hx10} showed that the problem determining bondage
number for general graphs is NP-hard. However, the bondage number
has received considerable attention in~\cite{cd06, dhtv98, fjkr90,
frv03, hl94, hx06, hx07a, hx08, hx10, ksk05, ky00, ts97, zlm09}.
Much work focused on the bounds of the bondage number as well as the
restraints on particular classes of graphs. In particular, Fink {\it
et al.}~\cite{fjkr90} showed
$b(K_n)=\left\lceil\frac{n}{2}\right\rceil$ for an $(n-1)$-regular
graph $K_n$ of order $n\geqslant2$, $b(G)=n-1$ for  an
$(n-2)$-regular graph $G$ of order $n\geqslant2$, where $G$ is a
$t$-partite graph $K_{n_1,n_2,\ldots,n_t}$ with $n_1=n_2=\cdots =
n_t=2$ and $t=\frac{n}{2}$ for an even integer $n\ge 4$.

In this paper, we show that $b(G)=n-3$ for every $(n-3)$-regular graph $G$
of order $n\ge 4$.

\section{Main results}

\begin{lem}\label{lem2.1}
$\gamma(G)=2$ for any $(n-3)$-regular graph $G$ of order $n\geq 4$.
\end{lem}

\begin{pf}
Let $G$ be an $(n-3)$-regular graph of order $n\geq 4$. It is clear
that $\gamma(G)\ge 2$ since there exists no such a vertex that
dominates all vertices in $G$. We only need to construct a
dominating set of $G$ with two vertices. Let $x$ be any vertex, and
let $y$ and $z$ be the only two vertices not adjacent to $x$ in $G$.

If $yz\in E(G)$, let $D=\{x,y\}$. If $yz\notin E(G)$, then there is
a vertex $w$ adjacent to both $y$ and $z$ in $G$ since $n\ge 4$. Let
$D=\{x,w\}$. Then $D$ is a dominating set of $G$. Thus $\gamma(G)\le
2$ and hence $\gamma(G)=2$.
\end{pf}

\begin{lem}\label{lem2.2}
$b(G)\leq n-3$ for any $(n-3)$-regular graph $G$ of order $n\geq 4$.
\end{lem}

\begin{pf}
For a vertex $x$ in $G$, let $G'=G-x$. Then $\gamma(G')\ge 2$ since
any vertex $y$ in $G'$ can not dominate all vertices in $G'$. Thus,
$\gamma(G-E_G(x))\ge 3>2=\gamma(G)$ by Lemma~\ref{lem2.1}, which
implies that $b(G)\leq n-3$.
\end{pf}

\begin{lem}\label{lem2.3}
Let $G$ be an $(n-3)$-regular graph of order $n\ge 7$ and $B$ be a
minimum bondage set of $G$. If $|B|\le n-4$, then there are at most
two vertices $x$ and $y$ in $G$ such that $E_G(x)\cap B=E_G(y)\cap
B=\emptyset$.
\end{lem}

\begin{pf}
Let $G'=G-B$. Then $\gamma(G')=3$ since $B$ is minimum and
$\gamma(G)=2$ by Lemma~\ref{lem2.1}. Suppose to the contrary that
there are three vertices $x_1,x_2$ and $x_3$ such that
 \begin{equation}\label{e2.1}
 E_G(x_i)\cap B=\emptyset\ \ {\rm for\ each}\ i=1,2,3.
 \end{equation}
Let $X=\{x_1,x_2,x_3\}$. We claim that $X$ is a dominating set of
$G'$. In fact, if there is some vertex $u$ in $G'$ that is not
dominated by $X$, that is, $u$ is not adjacent to any vertex in $X$,
then $u$ has degree at most $(n-4)$ in $G$ by (\ref{e2.1}), a
contradiction.
Thus, since X is a minimum
dominating set of $G'$,there exists a vertex $y_{ij}\notin N_{G'}[x_i]\cup N_{G'}[x_j]$
in $G'$ for $1\le i\ne j\le 3$ since $X$ is a minimum dominating set
of $G'$. Let $Y=\{y_{12},y_{23},y_{13}\}$. Then $|Y|=3$ since the
vertex $x_i$ has degree $(n-3)$ in $G'$ by (\ref{e2.1}) for each
$i=1,2,3$.

Since $G$ is $(n-3)$-regular, $\{y_{12},y_{13}\}$,
$\{y_{12},y_{23}\}$, and $\{y_{23},y_{13}\}$ are the only two
vertices not adjacent to $x_1$, $x_2$ and $x_3$ in $G$,
respectively. Similarly, $\{x_1,x_2\}$, $\{x_2,x_3\}$ and
$\{x_3,x_1\}$ are the only two vertices not adjacent to $y_{12}$,
$y_{23}$ and $y_{13}$ in $G$, respectively.  Also since $G$ is
$(n-3)$-regular, $G[Y]=K_3$. Note that any vertex in $X$ can
dominate all vertices in $G$ except for some two vertices in $Y$. If
$G[Y]-B$ contains a vertex, say $y_{23}$, of degree two, then
$\{x_1,y_{23}\}$ is a dominating set of $G-B$, a contradiction.
Therefore, $G[Y]-B$ contains no vertices of degree two, which means
that
 \begin{equation}\label{e2.2}
 |E(G[Y])\cap B|\ge 2.
 \end{equation}

Let $S=V(G)\setminus (X\cup Y)$. Then $S\ne \emptyset$ since $n\ge
7$. By (\ref{e2.2}) and $|B|\le n-4$, we have that $|E_G(Y,S)\cap
B|\le n-6$. Then there is some $s\in S$ such that $|E_G(s,Y)\cap
B|\le 1$. Without loss of generality, assume that two edges
$sy_{12}$ and $sy_{13}$ are both not in $B$. Then $\{s,x_1\}$ is a
dominating set of $G'$ since $s$ can dominate two vertices $y_{12}$
and $y_{13}$ and $x_1$ can dominate all vertices except for
$y_{12}$ and $y_{13}$, which implies that $\gamma(G')\le 2$, a
contradiction.
\end{pf}

\begin{lem}\label{lem2.4}
Let $G$ be an {\rm ($n-3$)}-regular graph of order $n\ge 7$ and $B$
be a minimum bondage set of $G$. If there exists a vertex $x\in
V(G)$ such that $E_G(x)\cap B=\emptyset$, then $|B|=n-3$.
\end{lem}

\begin{pf}
By Lemma~\ref{lem2.2}, we only need to prove that $|B|\geq n-3$.
Let $x\in V(G)$ such that
$E_G(x)\cap B=\emptyset$, let $y,z$ be the only two vertices not adjacent
to $x$ in $G$, and $G'=G-B$. Then $\gamma(G')\ge 3$ by
Lemma~\ref{lem2.1}.

If $yz\notin E(G)$, then both $y$ and $z$ are adjacent to each vertex
$s\in V_1=V(G)\setminus \{x,y,z\}$ in $G$. Thus at least one edge in
$\{sy,sz\}$ belongs to $B$ otherwise $\{x,s\}$ is a dominating set
of $G'$. Then
 $$
 |B|\ge|E_G(\{y,z\},V_1)\cap B|\ge |V_1|=n-3.
 $$

If $yz\in E(G)$, then $yz\in B$ otherwise $\{x,y\}$ is a dominating
set of $G'$. Let $p$ and $q$ be the only vertex except $x$ not
adjacent to $y$ and $z$ in $G$, respectively. Then both $y$ and $z$
are adjacent to any vertex $t\in V_2=V(G)\setminus \{x,y,z,p,q\}$ in
$G$. Thus at least one of $\{ty,tz\}$ belongs to $B$  otherwise
$\{x,t\}$ is a dominating set of $G'$. Thus,
  \begin{equation}\label{e2.3}
 |E_G(\{y,z\},V_2)\cap B|\ge |V_2|=
 \left\{\begin{array}{ll}
 n-4\ & \ {\rm if}\ p=q;\\
 n-5\ & \ {\rm if}\ p\ne q.
 \end{array}\right.
 \end{equation}

If $p=q$ then, by (\ref{e2.3}), we have that
 $$
 |B|\ge |\{yz\}|+|E_G(\{y,z\},V_2)\cap B|\ge n-3.
 $$

If $p\ne q$, then $|(E_G(p)\cup E_G(q))\cap B|\ge 1$ since
$E_G(x)\cap B=\emptyset$ and there are at most two vertices $u,v\in
V(G)$ such that $E_G(u)\cap B=E_G(v)\cap B=\emptyset$ by
Lemma~\ref{lem2.3}. Thus, by (\ref{e2.3}), we have that
 $$
 |B|\ge |\{yz\}|+|E_G(\{y,z\},V_2)\cap B|+|(E_G(p)\cup E_G(q))\cap
 B|\ge n-3.
 $$
The lemma follows.\end{pf}

\begin{thm}\label{thm2.1}
 $b(G)=n-3$ for any {\rm ($n-3$)}-regular graph $G$ of order $n\ge 4$.
\end{thm}

\begin{pf}
We first consider $n\in\{4,5,6\}$. If $n=4$, then $G=K_2+K_2$, so
$b(G)=1$. If $n=5$, then $G=C_5$, thus $b(G)=2$. Assume $n=6$ below.

Let $x$ be a vertex and $y,z$ be the only two vertices not adjacent to $x$
in $G$. It is easy to verify that
 $$
 G=\left\{\begin{array}{ll}
 C_3\times K_2\ &\ {\rm if}\ yz\in E(G);\\
 K_{3,3}\ &\ {\rm if}\ yz\notin E(G),
 \end{array}\right.
 $$
and so $b(G)=3$.

We now assume $n\ge 7$ in the following discussion.

Let $B$ be a minimum bondage set of $G$ and $G'=G-B$. Then
$\gamma(G')=3$ by Lemma~\ref{lem2.1}. If there exists a vertex $x\in
V(G)$ such that $E_G(x)\cap B=\emptyset$, then $|B|=n-3$ by
Lemma~\ref{lem2.4}. We now assume that $E_G(x)\cap B\ne \emptyset$
for every vertex  $x\in V(G)$. By Lemma~\ref{lem2.2}, $|B|\leq n-3$.
Next, we prove that $|B|\ge n-3$. Then there exists a vertex $x\in
V(G)$ such that $|E_G(x)\cap B|=1$.

Let $xw\in B$, $y$ and $z$ be the only two vertices not adjacent to
$x$ in $G$. Let $p$ and $q$ be the only two vertices not adjacent to
$w$ in $G$. We claim that for any vertex $x'\in V(G)\setminus
\{x,y,z,w\}$,
  \begin{equation}\label{e2.4}
 |E_G(\{w,y,z\},x')\cap B|\ge 1\ \ {\rm if}\ \{wx',yx',zx'\}\subseteq
 E(G).
 \end{equation}
To see this, note that if $|E_G(\{w,y,z\},x')\cap B|=\emptyset$, then $\{x,x'\}$
is a dominating set of $G'$ since $w,y$ and $z$ can be dominated by
$x'$ and others can be dominated by $x$ in $G'$, a contradiction.

We now prove that $|B|\ge n-3$ by considering the following three
cases.

\begin{description}

\item [Case 1]  $\{y,z\}=\{p,q\}$.

In this case, $yz\in E(G)$ and $x'$ is adjacent to every vertex in
$\{w,y,z\}$ for any vertex $x'\in V_1=V(G)\setminus \{x,y,z,w\}$ in
$G$. By (\ref{e2.4}), $|E_G(\{w,y,z\},x')\cap B|\ge 1$, and so
$|E_G(\{w,y,z\},V_1)\cap B|\geq |V_1|=n-4$. Thus
$$
 \begin{array}{rl}
 |B|&\geq |\{xw\}|+|E_G(\{w,y,z\},V_1)\cap B|\\
           &\geq n-3.
  \end{array}
$$

\item [Case 2]  $|\{y,z\}\cap \{p,q\}|=1$. Without loss of generality, let $p=y$.

In this case, $yz,wz\in E(G)$ and hence $|E_G(z,\{y,w\})\cap B|\geq
1$, for otherwise $\{x,z\}$ is a dominating set of $G'$ since $\{y,w\}$
can be dominated by $z$ and others can be dominated by $x$ in $G'$.
Let $r$ be the only vertex except $x$ not adjacent to $z$ in $G$.
Thus, $x'$ is adjacent to every vertex in $\{w,y,z\}$ for any vertex
$x'\in V_2=V(G)\setminus \{x,y,z,w,q,r\}$ in $G$. By (\ref{e2.4}),
$|E_G(\{w,y,z\},x')\cap B|\ge 1$, and so $|E_G(\{w,y,z\},V_2)\cap
B|\geq |V_2|=n-6$. Thus
 $$
 \begin{array}{rl}
 |B|&\geq |\{xw\}|+|E_G(\{w,y,z\},V_2)\cap B|\\
            &+|E_G(z,\{y,w\})\cap B|+|(E_G(q)\cup E_G(r))\cap B|\\
           &\geq n-3.
  \end{array}
$$

\item [Case 3]  $\{y,z\}\cap \{p,q\}=\emptyset$.

In this case, $wy,wz\in E(G)$.

\begin{description}

\item [subcase 3.1] $yz\notin E(G)$.

In this case, $|E_G(w,\{y,z\})\cap B|\ge 1$, for otherwise $\{x,w\}$ is a
dominating set of $G'$. Note that $x'$ is adjacent to every vertex
in $\{w,y,z\}$ for any vertex $x'\in V_3=V(G)\setminus
\{x,y,z,w,p,q\}$ in $G$. By (\ref{e2.4}), $|E_G(\{w,y,z\},x')\cap
B|\ge 1$, and so $|E_G(\{w,y,z\},V_2)\cap B|\geq |V_2|=n-6$. Thus
$$
 \begin{array}{rl}
  |B|&\geq |\{xw\}|+|E_G(\{w,y,z\},V_3)\cap B|\\
            &+|E_G(w,\{y,z\})\cap B|+|(E_G(p)\cup E_G(q))\cap B|\\
           &\geq n-3.
  \end{array}
$$

\item [Subcase 3.2]  $yz\in E(G)$.

In this case, $|E(G[\{w,y,z\}])\cap B|\ge 2$, since otherwise $x$ and one
vertex in $\{w,y,z\}$ consist of a dominating set of $G'$. Let $r$
and $s$ be the only vertex except $x$ not adjacent to $y$ and $z$ in
$G$, respectively. Note that $x'$ is adjacent to every vertex in
$\{w,y,z\}$ for any vertex $x'\in V_4=V(G)\setminus
\{x,y,z,w,p,q,r,s\}$ in $G$. By (\ref{e2.4}),
$|E_G(\{w,y,z\},x')\cap B|\ge 1$, and so
 $$
 |E_G(\{w,y,z\},V_4)\cap B|\geq
 \left\{\begin{array}{ll}
 n-6\ &\ {\rm if}\ |\{r,s\}\cup \{p,q\}|=2;\\
 n-7\ &\ {\rm if}\ |\{r,s\}\cup \{p,q\}|=3;\\
 n-8\ &\ {\rm if}\ |\{r,s\}\cup \{p,q\}|=4.
 \end{array}\right.
 $$
Thus
$$
 \begin{array}{rl}
 |B|&\geq |\{xw\}|+|E_G(\{w,y,z\},V_4)\cap B|\\
 &+|E(G[\{w,y,z\}])\cap B|\\
  &+|(E_G(p)\cup E_G(q)\cup E_G(r)\cup E_G(s))\cap B|\\
           &\geq n-3.
  \end{array}
$$
\end{description}

\end{description}
The theorem follows.
\end{pf}


\begin{thebibliography}{99}

\bibitem{cd06}
K. Carlson and M. Develin, On the bondage number of planar and
directed graphs. {\it Discrete Mathematics}, {\bf 306} (8-9) (2006),
820-826.

\bibitem{dhtv98}
J. E. Dunbar, T. W. Haynes, U. Teschner, L. Volkmann, Bondage,
insensitivity, and reinforcement. {\it Domination in Graphs:
Advanced Topics} (T.W. Haynes, S.T. Hedetniemi, P.J. Slater eds.),
Marcel Dekker, New York, 1998, pp. 471-489.

\bibitem{fjkr90}
J. F. Fink, M. S. Jacobson, L. F. Kinch, J. Roberts, The bondage
number of a graph. {\it Discrete Mathematics}, {\bf 86} (1990),
47-57.

\bibitem{frv03}
M. Fischermann, D. Rautenbach and L. Volkmann, Remarks on the
bondage number of planar graphs. {\it Discrete Mathematics}, {\bf
260} (2003), 57-67.

\bibitem{hl94}
B. L. Hartnell and D. F. Rall, Bounds on the bondage number of a
graph. {\it Discrete Mathematics}, {\bf 128} (1994), 173-177.





\bibitem{hx10}
F.-T. Hu and J.-M. Xu, Complexity of bondage and reinforcement. A
manuscript submitted to {\it Discrete Applied Mathematics}, 2010.

\bibitem{hx06}
J. Huang and J.-M. Xu, The bondage numbers of extended de Bruijn and
Kautz digraphs. {\it Computer and Mathematics with Applications},
{\bf 51} (6-7) (2006), 1137-1147.

\bibitem{hx07a}
J. Huang and J.-M. Xu, The bondage number of graphs with small
crossing number. {\it Discrete Mathematics}, {\bf 307} (15) (2007),
1881-1897.

\bibitem{hx08}
J. Huang and J.-M. Xu, The bondage numbers and efficient dominations
of vertex-transitive graphs. {\it Discrete Mathematics}, {\bf 308}
(4) (2008), 571-582.

\bibitem{hx10}
J. Huang and J.-M. Xu, Domination and total dominatoin contraction
numbers of graphs. {\it Ars Combinatoria}, {\bf 94} (2010), 431-443.


\bibitem{ksk05}
L.-Y Kang, M. Y. Sohn and H. K. Kim, Bondage number of the discrete
torus $C_n\times C_4$. {\it Discrete Mathematics}, {\bf 303} (2005),
80-86.

\bibitem{ky00}
L.-Y. Kang, J.-J. Yuan, Bondage number of planar graphs. {\it
Discrete Mathematics}, {\bf 222} (2000), 191-198.

\bibitem{ts97}
U. Teschner, New results about the bondage number of a graph. {\it
Discrete Math}. {\bf 171} (1997), 249-259.




\bibitem{x03}
J.-M. Xu, {\it Theory and Application of Graphs}. Kluwer Academic
Publishers, Dordrecht/Boston/London, 2003.

\bibitem{zlm09}
X. Zhang, J. Liuy and J.-X. Meng, The bondage number in complete
$t$-partite digraphs. {\it Information Processing Letters}, {\bf
109}(16) (2009), 997-1000.

\end{thebibliography}
\end{document}